\documentclass[aop,seceqn,MSNbibl,citesort,dvips]{arximspdf}
\usepackage{mathrsfs}

\doi{10.1214/09-AOP519}
\volume{38}
\issue{4}
\pubyear{2010}
\firstpage{1570}
\lastpage{1582}

\makeatletter

\newproclaim{remark}{Remark}[section]
\newtheorem{lemma}[remark]{Lemma}
\newtheorem{proposition}[remark]{Proposition}
\newtheorem{corollary}[remark]{Corollary}
\newtheorem{theorem}[remark]{Theorem}
\def\bsuffix #1{#1}
\makeatother

\begin{document}
\begin{frontmatter}

\title{Volume growth and escape rate of Brownian motion on a complete
Riemannian manifold}
\runtitle{Volume growth and escape rate of Brownian motion}

\begin{aug}
\author[A]{\fnms{Elton P.} \snm{Hsu}\corref{}\ead[label=e1]{ehsu@math.northwestern.edu}} and
\author[B]{\fnms{Guangnan} \snm{Qin}\ead[label=e2]{gn\_qin@amss.ac.cn}}
\runauthor{E. P. Hsu and G. Qin}
\affiliation{Northwestern University and Chinese Academy of Sciences}
\address[A]{Department of Mathematics\\
Northwestern University\\
Evanston, Illinois 60208\\
USA\\
\printead{e1}}
\address[B]{Institute of Applied Mathematics\\
Academy of Mathematics\\ \quad and Systems Science\\
Chinese Academy of Sciences\\
Beijing\\
China\\
\printead{e2}}
\end{aug}

\received{\smonth{10} \syear{2009}}
\revised{\smonth{11} \syear{2009}}

%
\begin{abstract}
We give an effective upper escape rate function for Brownian motion on
a complete Riemannian manifold in terms of the volume growth of the
manifold. An important step in the work is estimating the small tail
probability of the crossing time between two concentric geodesic
spheres by reflecting Brownian motions on the larger geodesic ball.
\end{abstract}

%
\setattribute{keyword}{AMS}{AMS 2000 subject classification.}

\end{frontmatter}

\section{Introduction}\label{sec1}

Let $M$ be a Riemannian manifold and $p_M (t,x,y)$ the (minimal) heat
kernel on $M$. By definition, the latter is the fundamental solution of
the heat operator
\[
{\mathscr L}_M = \frac{\partial}{\partial t} - \frac12\Delta_M,
\]
where $\Delta_M$ is the Laplace--Beltrami operator on $M$. A Riemannian
manifold is stochastically complete if
\[
\int_M p_M (t,x,y) \,dy = 1
\]
for some, hence all, $(x,t)\in M\times(0,\infty)$. In other words, $M$
is stochastically complete if the heat kernel is conserved. Let
${\mathbb P}_x$
be the law of Brownian motion on $M$ starting from $x$ and let $e$ be
the lifetime (or explosion time) of the Brownian motion. We then have
\[
{\mathbb P}_x\{ e> t\}= \int_M p (t,x,y) \,dy.
\]
Therefore, $M$ is stochastically complete if and only if
%
\begin{equation}\label{notexplode}
{\mathbb P}_x\{ e = \infty\}= 1,
\end{equation}
that is, Brownian motion on $M$ does not explode. Finding geometric
conditions for stochastic completeness is an old geometric problem. The
problem has been attacked using both analytic and probabilistic
methods. Early works on this problem (Karp and Li
\cite{KarpL}, Yau \cite{Yau}, Hsu \cite{Hsu1} and Varopoulos \cite{Varopoulos}) impose
lower bounds on the Ricci curvature. In particular, in the last two
works, it was shown that if there is a strictly positive function
$\kappa(r)$ such that the Ricci curvature of $M$ on the geodesic ball
$B(r)$ is bounded from below by $-\kappa(r)$ and
\[
\int_0^\infty\frac{dr}{\sqrt{\kappa(r)}} = \infty,
\]
then $M$ is stochastically complete. In 1986, Grigor'yan
\cite{Grigoryan} found the following sufficient condition for stochastic
completeness solely in terms of the volume growth function of the manifold:
%
\begin{equation}\label{itest}
\int_1^\infty\frac{r \,dr}{{\ln}\vert B(r)\vert} = \infty.
\end{equation}

According to (\ref{notexplode}), a Riemannian manifold $M$ is
stochastically complete if Brownian motion does not escape to infinity
(in the one-point compactification) in finite time. Traditionally in
probability theory (see, e.g., It\^o and McKean \cite{ItoM} and Shiga
and Watanabe \cite{ShigaW}), one often looks for upper functions for
the escape rate of a diffusion process. The classical Khinchin law of
iterated logarithm for one-dimensional standard Brownian motion is the
most celebrated case. More generally, let $r_t = d(X_t, x)$ be the
radial process of Brownian motion $X$ on $M$. An increasing function
$\psi(t)$ is called an upper rate function if
\[
{\mathbb P}_x\{ r_t\le\psi (t)\mbox{ for all sufficiently large }
t\}= 1.
\]
Finding an upper rate function is a more refined problem than proving
stochastic completeness since the existence of an upper rate function
implies stochastic completeness. Various explicit upper rate function
for Brownian motion on a complete Riemannian manifold have been
obtained under concrete volume growth assumptions (see Grigor'yan \cite
{Grigoryan2}, Grigor'yan and Kelbert \cite{GrigoryanK} and
Takeda \cite{Takeda,Takeda1}). More recently, Grigor'yan and Hsu \cite
{GrigoryanH} showed that the inverse function of the increasing function
\[
\phi_1(R) = \int_1^R\frac{r \,dr}{{\ln}\vert B(r)\vert}
\]
related to the integral test (\ref{itest}) for stochastic completeness is
essentially an upper rate function for Brownian motion on $M$. While
this result gives an upper rate function of a very general form, it was
proven under the additional geometric assumption that $M$ is a
Cartan--Hadamard manifold, that is, a simply connected, geodesically
complete Riemannian manifold with nonpositive sectional curvature.

The main purpose of the present work is to obtain an escape rate
function based solely on the volume growth of the underlying manifold.
We introduce the following increasing function:
\[
\phi(R) = \int_6^R\frac{r \,dr}{{\ln}\vert B(r)\vert+ \ln\ln r}.
\]
We will show in Theorem \ref{main} that under the sole assumption that
$M$ is a complete Riemannian manifold, the inverse function of $\phi$
is essentially an upper rate function for Brownian motion on $M$.
\begin{remark} \label{nobetter} The difference between the functions
$\phi_1(r)$ and $\phi(r)$ is that we have introduced an extra term
$\ln\ln r$ in the latter function. This addition, resulting from an
attempt to remove the extraneous geometric condition, is fully
justified on several grounds. First, the integral test (\ref{itest}) for
stochastic completeness, which can be written as $\phi_1(\infty) =
\infty$, is equivalent to the condition $\phi(\infty) = \infty$.
Second, our new upper rate function implies all explicit upper rate
functions existent in the literature to date (see Corollaries
\ref{oldrates} and \ref{newrate}). Third, in general, the radial
process of a Brownian motion on $M$ has the form
\[
r_t = \beta_t + \frac12\int_0^t \Delta_M r(X_t) \,dt - L_t,
\]
where $\beta$ is a standard one-dimensional Brownian motion and $L$ is
a local time on the cut locus $C(x)$ of the point $x$ (see
Hsu \cite{Hsubook} and Kendall \cite{Kendall}). In the absence of any further
geometric assumptions, we do not expect to obtain an upper rate
function (up to a multiplicative constant) for the process $r_t$ better
than the upper rate function of a standard Brownian motion $\psi(t) =
C\sqrt{t\ln\ln t}$. This rate function cannot be achieved without the
presence of the additional term $\ln\ln r$ in the function $\phi$.
\end{remark}

Our method has two key steps. In the first, we follow Hsu \cite{Hsu1}
and Grigor'yan and Hsu \cite{GrigoryanH} and, by using the
Borel--Cantelli lemma, reduce the problem of seeking an upper rate
function to the problem of estimating the small tail probability of the
crossing time between two concentric geodesic spheres (Lemma \ref
{usebc}). In the second key step, instead of estimating the small tail
probability by using an analytic approach, as in Grigor'yan and
Hsu \cite{GrigoryanH}, under the assumption that the manifold in
question is Cartan--Hadamard, we modify the way
the Lyons--Zheng decomposition for
reflecting Brownian motion is used in Takeda
\cite{Takeda,Takeda1}. The volume $\vert B(r)\vert$ of the
geodesic ball $B(r)$ appears naturally in this step because the uniform
distribution on the ball with respect to the Riemannian volume measure
is the invariant measure of reflecting Brownian motion on $B(r)$. The
additional term mentioned above, $\ln\ln r$, is a consequence of
dealing with the Brownian motion adapted to the time-reversed
filtration in the Lyons--Zheng decomposition.

\section{Basic estimates on crossing times}

Let $M$ be a geodesically complete Riemannian manifold and ${\mathscr P}(M)$
the path space over $M$. Let $X$ be the canonical coordinate process on
the path space ${\mathscr P}(M)$ over $M$, that is, $X_t(\omega) =
\omega_t$
for $\omega\in{\mathscr P}(M)$. If $x\in M$, then we use ${\mathbb
P}_x$ to denote the
law of Brownian motion on $M$ starting from $x$. The radial process is
$r_t = d(X_t, x)$, that is, the Riemannian distance from $x$ to $X_t$,
the position of Brownian motion at time $t$. A nonnegative increasing
function $R\dvtx {\mathbb R}_+\rightarrow{\mathbb R}_+$ is called an
upper rate function
for Brownian motion on $M$ if
\[
\mathbb{P}_x\{ r_t \leq R( t) \mbox{ for all sufficiently large
}t\}=1.
\]

Let $\{ R_n\}$ be a strictly increasing sequence of positive
numbers, to be chosen later, and define a sequence of stopping times as follows:
\[
\tau_n = \inf\{ t\dvtx r_t = R_n\}.
\]
(We adopt the convention that $\inf\varnothing= \infty$.) Each of these
is the first time the Brownian motion $X$ reaches the corresponding
geodesic sphere,
\[
S(R_n) = \{ x\in M\dvtx d(x, o) = R_n\}.
\]
The difference $\tau_n - \tau_{n-1}$, if well defined, is the amount of
time Brownian motion takes to cross from $S(R_{n-1})$ to $S(R_n)$. The
basic idea of Grigor'yan and Hsu \cite{GrigoryanH} (see also Hsu
\cite{Hsu1}) for controlling the rate of escape of Brownian motion is to
give a good upper bound for the small tail probability ${\mathbb P}_x\{
\tau
_n-\tau_{n-1}\le t_n\}$ for a suitably chosen sequence $\{ t_n\}$
of time steps. If the sum of these probabilities converges, then the
Borel--Cantelli lemma shows that for sufficiently large $n$, Brownian
motion $X$ has to wait roughly until at least
\[
T_n = \sum_{k=1}^n t_k
\]
to reach the sphere $S(R_n)$ or, equivalently, $r_t\le R_n$ for all
$t\le T_n$. This, after some technical manipulations (see Section
\ref{getrate}), will give an upper escape rate function.

We now use the idea of Takeda \cite{Takeda,Takeda1} to estimate
the small tail probability ${\mathbb P}_x\{\tau_n -\tau_{n-1}\le t_n\}
$ by
using the Lyons--Zheng decomposition \cite{LyonsZ} of reflecting
Brownian motion starting from the uniform distribution on a geodesic
ball. For an open set $B\subset M$, we denote by ${\mathbb P}_B$ the
law of
Brownian motion starting from the uniform distribution on $B$, that is,
\[
{\mathbb P}_B = \frac1{\vert B\vert}\int_B{\mathbb P}_x \,dx.
\]
Likewise, we use ${\mathbb Q}_B$ to denote the law of reflecting Brownian
motion on $B$ starting from the same uniform distribution. Let $B_n =
B(R_n)$ be the geodesic ball of radius $R_n$ centered at $x$. In order
to take advantage of the volume growth condition, we consider
the probability ${\mathbb P}_{B_1}\{\tau_n-\tau_{n-1}\le t_n\}$
instead of
${\mathbb P}_x\{\tau_n-\tau_{n-1}\le t_n\}$.
Recall that $\tau_n$ is the first time the process $X$ reaches the
boundary $S(R_n)$ of the geodesic ball $B(R_n)$. Before reaching the
boundary, Brownian motion and reflecting Brownian motion have the same
law. Therefore, if $C\in{\mathscr B}_{\tau_n}$ is an event which is measurable
up to time $\tau_n$, then ${\mathbb P}_{B_n}(C) = {\mathbb
Q}_{B_n}(C)$. From
\[
{\mathbb P}_{B_1}(C) = \frac1{\vert B_1\vert}\int_{B_1}{\mathbb
P}_x(C) \,dx\le
\frac1{\vert B_1\vert}\int_{B_n}{\mathbb P}_x(C) \,dx=\frac{\vert
B_n\vert
}{\vert B_1\vert} {\mathbb P}_{B_n}(C),
\]
we have
%
\begin{equation}\label{bmtorbm}
{\mathbb P}_{B_1}(C) \le\frac{\vert B_n\vert}{\vert B_1\vert}
{\mathbb Q}
_{B_n}(C),\qquad C\in{\mathscr B}_{\tau_n}.
\end{equation}
We apply this inequality to the event
%
\begin{equation}\label{notc}
C = \{\tau_n-\tau_{n-1}\le t_n\}.
\end{equation}

Now, according to the Lyons--Zheng decomposition \cite{LyonsZ}, on a
fixed time horizon $[0, T_n]$, the radial process can be decomposed as
the difference
%
\begin{equation}\label{lyonszheng}
r_t- r_0 = \frac{B_t}2 -\frac{\tilde B_{T_n} - \tilde
B_t}2,
\end{equation}
where $B$ is a standard Brownian motion adapted to the natural
filtration $msb_*={\mathscr B}({\mathscr P}(M))_*$ of the path space
${\mathscr P}(M)$ and
$\tilde B$ is also a standard Brownian motion, but adapted to the
reversed filtration $\tilde{\mathscr B}_*$ defined by
\[
\tilde{\mathscr B}_t = \sigma\{ X_{T_n-s}\dvtx 0\le s\le t\},\qquad 0\le
t\le T_n.
\]
The advantage of such a decomposition is obvious, for we have
eliminated from consideration the bounded variation component of the
radial process, which can be rather complicated. The price is that we
have to deal with a Brownian motion not adapted to the original
filtration. Another complication is that the decomposition cannot be
applied directly to the event (\ref{notc}) because it may go beyond the
fixed time horizon $[0, T_n]$. In order to remedy this situation, we
will use a slightly modified event
\[
C_n = \{\tau_n - \tau_{n-1}\le t_n, \tau_n\le T_n\}
\]
in the inequality (\ref{bmtorbm}). Fortunately, this additional
restriction $\{\tau_n\le T_n\}$ will not be an obstacle for us, as
shown by the following simple observation.
\begin{lemma} \label{usebc} Let $C_n$ be defined as above. Suppose that
$\sum_{n=1}^\infty{\mathbb P}(C_n)<\infty$. Then, with probability
1, there
exists $T_{-1}$ such that $\tau_n\ge T_n - T_{-1}$ for all $n$.
\end{lemma}
\begin{pf}
By the Borel--Cantelli lemma, the probability that the events $\{
C_n\}$ happen infinitely often is 0. Therefore, with probability 1,
there exists $n_0$ such that for all $n\ge n_0$, either $\tau_n-\tau
_{n-1}\ge t_n$ or $\tau_n\ge T_n$. We show, by induction, that $\tau
_n\ge T_n - T_{n_0}$ holds for all $n$. If $1\le n\le n_0$, then $\tau
_n\ge0\ge T_n - T_{n_0}$. Suppose that $\tau_n\ge T_n - T_0$ for an
$n\ge n_0$. If $\tau_{n+1}\ge T_{n+1}$, then, trivially, $\tau
_{n+1}\ge
T_{n+1}- T_{n_0}$. Otherwise, $\tau_{n+1}-\tau_n\ge t_{n+1}$ and
\[
\tau_{n+1} = \tau_{n+1} - \tau_n +\tau_n \ge t_{n+1} +T_n - T_{n_0} =
T_{n+1} -T_{n_0}.
\]
This completes the proof.
\end{pf}

We now prove the main estimate for the crossing time $\tau_n - \tau_{n-1}$.
\begin{proposition} Let $\tau_n$ be the first hitting time of the
sphere $S(R_n)$ and $r_n = R_n - R_{n-1}$. There then exists a constant
$C$ such that
\[
{\mathbb P}_{B_1}\{\tau_n - \tau_{n-1}\le t_n, \tau_n\le T_n\}\le
\frac
{\vert B_n\vert}{\vert B_1\vert}\frac{C}{\sqrt{\pi t_n}}\frac{T_n}{r_n}
e^{-r_n^2/8t_n}.
\]
\end{proposition}
\begin{pf} The event $\{\tau_n - \tau_{n-1}\le t_n\}$ implies
the event
%
\begin{equation}\label{oscillation}
\Bigl\{\sup_{0\le s\le t_n }(r_{\tau_{n-1}+s} - r_{\tau_{n-1}})\ge r_n\Bigr\}.
\end{equation}
Now, from the decomposition (\ref{lyonszheng}), we have
\[
2(r_{\tau_{n-1}+s}-r_{\tau_{n-1}}) = B_{\tau_{n-1}+s} - B_{\tau_{n-1}}+
\tilde B_{\tau_{n-1}+s} - \tilde B_{\tau_{n-1}}.
\]
Since $\tau_{n-1}$ is a stopping time with respect to the natural
filtration ${\mathscr B}_*$, the first term on the right-hand side is a
Brownian motion in time $s$ starting from 0. This is not so for the
second term because $\tau_{n-1}$ is not a stopping time with respect to
the filtration $\tilde{\mathscr B}_*$ of the reversed process. However,
for any $s\le t_n$ such that $\tau_{n-1}\le T_n$, taking $k$ such that
$(k-1)t_n\le\tau_{n-1}\le kt_n$, we see that both $\tau_{n-1}$ and
$\tau
_{n-1}+s$ lie in the interval $[(k-1)t_n, (k+1)t_n]$.
From
\[
r_{\tau_{n-1}+s} - r_{\tau_{n-1}} = r_{\tau_{n-1}+s} - r_{kt_n} +
r_{kt_n} -r_{\tau_{n-1}},
\]
the event (\ref{oscillation}) is contained in the union of the
$[T_n/t_n]+1$ events
\[
\biggl\{{\sup_{\vert s\vert\le t_n}}\vert r_{kt_n+s}-r_{kt_n}\vert\ge\frac
{r_n}2\biggr\},\qquad 1\le k\le\biggl[\frac{T_n}{t_n} \biggr]+1.
\]
Using
\[
r_{kt_n +s}-r_{kt_n} = \frac{B_{kt_n +s} - B_{kt_n}}2 + \frac{\tilde
B_{kt_n +s} - \tilde B_{kt_n}}2,
\]
we see that the event $\{\tau_n - \tau_{n-1}\le t_n, \tau_n\le
T_n\}$ is also contained in the union of the following $2[T_n/t_n]+2$ events:
\[
\biggl\{{\sup_{\vert s\vert\le t_n}} | B_{kt_n+s} - B_{kt_n} |\ge
\frac{r_n}2\biggr\}
\]
and
\[
\biggl\{{\sup_{\vert s\vert\le t_n}} |\tilde B_{kt_n+s} - \tilde
B_{kt_n} |\ge\frac{r_n}2\biggr\}
\]
for $1\le k\le[T_n/t_n]+1$. Under the probability ${\mathbb Q}_{B_n}$, these
events have the same probability,
\[
{\mathbb P}\biggl\{{\sup_{\vert s\vert\le t_n}}\vert B_{t_n+s}-B_{t_n}\vert
\ge
\frac
{r_n}2\biggr\}\le2{\mathbb P}\biggl\{{\sup_{0\le s\le t_n}}\vert B_s\vert\ge\frac
{r_n}2\biggr\}\le\frac{C\sqrt t_n}{r_n} e^{-r_n^2/8t_n}.
\]
The probability ${\mathbb Q}_{B_n}\{\tau_n - \tau_{n-1}\le t_n,
\tau_n\le
T_n\}$ is bounded from above by $2[T_n/t_n]+2\le4T_n/t_n$ times the
above probability. The desired inequality now follows immediately from
this and the inequality [see (\ref{bmtorbm})]
\[
{\mathbb P}_{B_1}\{\tau_n -\tau_{n-1}\le t_n, \tau_n\le T_n\}\le
\frac
{\vert B_n\vert}{\vert B_1\vert} {\mathbb Q}_{B_n}\{\tau_n-\tau
_{n-1}\le
t_n, \tau_n\le T_n\}.
\]
\upqed\end{pf}

\section{Total crossing time}\label{total}

In the preceding section, we have found an upper bound for the
probability ${\mathbb P}_{B_1}\{\tau_n-\tau_{n-1}\le t_n, \tau
_n\le T_n\}
$. We are still free to choose the upper bounds $t_n$ of the crossing
times $\tau_n - \tau_{n-1}$ and the radii $R_n$ of the expanding
geodesic balls $B(R_n)$. We need to choose them so that the series
\begin{eqnarray*}
&&\sum_{n=1}^\infty{\mathbb P}_{B_1}\{\tau_n -\tau_{n-1}\le t_n,
\tau
_n\le
T_n\}\\
&&\qquad \le\frac C{\vert B_1\vert}\sum_{n=1}^\infty\frac
{T_n}{\sqrt
{t_n}r_n}\exp\biggl[{\ln}\vert B(R_n)\vert- \frac{r_n^2}{8t_n}
\biggr]\\
&&\qquad\le\frac{C_1}{\vert B_1\vert}\sum_{n=1}^\infty\frac
{T_n}{r^2_n}\exp
\biggl[{\ln}\vert B(R_n)\vert- \frac{r_n^2}{16t_n} \biggr]
\end{eqnarray*}
converges and the Borel--Cantelli lemma can be applied. The obvious
choice is for $t_n$ equal to a small multiple of $r_n^2/{\ln}\vert
B(R_n)\vert$, as was adopted in Grigor'yan and Hsu \cite{GrigoryanH}.
However, this choice will not enable us to eliminate the extra factor
$T_n/r_n^2$, whose presence can be traced back to the Brownian motion
$\tilde B$ adapted to the reverse filtration $\tilde{\mathscr B}_*$ in the
Lyons--Zheng decomposition (\ref{lyonszheng}). We diminish the obvious
choice by letting
%
\begin{equation}\label{timeincrement}
t_n = \frac1{32} \frac{r_n^2}{{\ln}\vert B(R_n)\vert+ h(R_n)}
\end{equation}
with a strictly increasing function $h$ to be determined. If we assume,
without loss of generality, that $B(R_1)\ge1$ and $h(R_1)\ge1$, then
$t_n\le r^2_n/32$. If we further assume that the sequence $\{ r_n\}
$ is increasing, then there is an obvious bound
\[
32T_n\le\sum_{k = 1}^n r^2_k\le\sum_{k=1}^n r_kr_n=R_nr_n.
\]
It follows that
\[
{\mathbb P}_{B_1}\{\tau_n-\tau_{n-1}\le t_n, \tau_n\le T_n\}\le
\frac
{C_2}{\vert B_1\vert^2} \frac{R_n}{r_n} e^{-2h(R_n)}.
\]
It remains to choose the radii $R_n$ and the function $h$ such that
%
\begin{equation}\label{divcon}
\sum_{n=1}^\infty t_n = \infty\quad\mbox{and}\quad \sum
_{n=1}^\infty
\frac{R_n}{r_n} e^{-2h(R_n)}<\infty
\end{equation}
under the integral condition
%
\begin{equation}\label{itest1}
\int_1^\infty\frac{r \,dr}{{\ln}\vert B(r)\vert} = \infty.
\end{equation}
The divergence of the above integral is to be linked to the divergence
of the total crossing time in (\ref{divcon}). This leads to the natural
requirement that
\[
r_n^2\ge CR_n(R_{n+1}- R_n).
\]
This requirement can be fulfilled by setting $R_n = 2^n$ with $C =
1/4$. From (\ref{timeincrement}), we have
%
\begin{equation}\label{accumulate}
T_n = \sum_{k =1}^n t_k = \frac1{128}\sum_{k=1}^n \frac
{R_k(R_{k+1}-R_k)}{{\ln}\vert B(R_k)\vert+h(R_k)}\ge\frac1{256}\int
_{R_1}^{R_{n+1}}\frac{r \,dr}{{\ln}\vert B(r)\vert+ h(r)},\hspace*{-28pt}
\end{equation}
which seems to fall slightly short of the condition (\ref{itest1}). The
apparently disadvantageous situation can be salvaged by first
looking at a typical candidate for the function $h$. From our choice of
$R_n = 2^n$, we have $R_n/r_n = 2$ and the convergence of the total
probability in (\ref{divcon}) becomes
\[
\sum_{n=1}^\infty e^{-2h(2^n)}<\infty.
\]
This leads to the choice $h(R) = \ln\ln R$. We have the following
simple observation.
\begin{lemma} \label{alsodiverge} Let $f$ be a positive, nondecreasing
and continuous function on $[0,+\infty)$ such that
%
\begin{equation}
\int^{\infty}_3\frac{r \,dr}{f(r)}=\infty.
\end{equation}
Then
\[
\int^{\infty}_3\frac{r \,dr}{f(r)+\ln\ln r}=\infty.
\]
\end{lemma}
\begin{pf} Divide the integral into the sum of the integrals over the
intervals $[n-1, n]$ for $n\ge4$. Since $f$ is increasing, we have
\begin{eqnarray*}
\int^{\infty}_3\frac{r \,dr}{f(r)+\ln\ln r}&\ge&\sum_{n=4}^\infty
\frac
{n-1}{f(n)+\ln\ln n}\\
&\ge&\frac12\sum_{f(n)\ge\ln\ln n}
\frac{n-1}{f(n)}+\frac12\sum_{f(n)< \ln\ln n}
\frac{n-1}{\ln\ln n}.
\end{eqnarray*}
Since $(n-1)/\ln\ln n\ge1$ for all sufficiently large $n$, if the
second sum has infinitely many terms, then it is clearly diverges;
otherwise, $f(n)\ge\ln\ln n$ for all sufficiently large $n$ and we
have, for some $ n_0$,
\begin{eqnarray*}
\int^{\infty}_3\frac{r \,dr}{f(r)+\ln\ln r}&\ge&\frac12 \sum
_{n=n_0}^\infty\frac{n-1}{f(n)}\\
&\ge&\frac12\sum_{n=n_0}^\infty
\frac{n-1}{n+1}\int_n^{n+1}\frac{r \,dr}{f(r)}\\
&\ge&\frac12\frac{n_0-1}{n_0+1}\int_{n_0}^\infty\frac{r \,dr}{f(r)}.
\end{eqnarray*}
This completes the proof.
\end{pf}

We are now in a position to bound the range of Brownian motion on a
finite time interval.
\begin{proposition} \label{frange} Let $R_n = 2^n$ and
\[
T_n = \frac1{128}\sum_{k=1}^n \frac{R_k(R_{k+1}-R_k)}{{\ln}\vert
B(R_k)\vert+h(R_k)}.
\]
Then, with probability 1, there exists $T_{-1}$ such that
$\sup_{0\le t\le T_n-T_{-1}} r_t\le2^n$
for all~$n$.
\end{proposition}
\begin{pf} By our choice of $R_n$,
\[
\sum_{n=1}^\infty{\mathbb P}_{B_1}\{\tau_n-\tau_{n-1}\le t_n, \tau
_n\le
T_n\}<\infty.
\]
By the Borel--Cantelli lemma and Lemma \ref{usebc}, with probability 1,
there exists $T_{-1}$ such that $\tau_n\ge T_n- T_{-1}$. However,
$\tau
_n$ is the hitting time of the sphere $S(2^n)$, hence
$\sup_{0\le t\le T_n -T_{-1}} r_t\le2^n$ for sufficiently large $n$.
\end{pf}

An easy consequence of the above result is a probabilistic proof of
Grigor'yan's criterion for stochastic completeness.
\begin{corollary}[(Grigor'yan \cite{Grigoryan})] Suppose that $M$ is a
complete Riemannian manifold and $B(R)$ its geodesic ball of radius $R$
centered at a fixed point. If
\[
\int_1^\infty\frac{r \,dr}{{\ln}\vert B(r)\vert} = \infty,
\]
then $M$ is stochastically complete.
\end{corollary}
\begin{pf} By Lemma \ref{alsodiverge}, under the assumption, we have
\[
T_n\ge\frac1{256}\int_{R_1}^{R_{n+1}}\frac{r \,dr}{{\ln}\vert
B(r)\vert+
h(r)}\rightarrow\infty
\]
as $n\rightarrow\infty$. By the above proposition, $\sup_{t\le
T}r_t<\infty$ for all finite $T$. Hence, Brownian motion does not
explode and $M$ is stochastic complete.
\end{pf}

\section{Upper rate function}\label{getrate}

Proposition \ref{frange} allows us to obtain an upper rate function in
terms of the volume growth function $\vert B(r)\vert$, as was similarly
done in Grigor'yan and Hsu \cite{GrigoryanH}.

Let
\[
\phi(R) = \int_6^R\frac{r \,dr}{{\ln}\vert B(r)\vert+\ln\ln r}.
\]
From (\ref{accumulate}), we have $(1/256)\phi(2^{n+1})\le T_n$.
Proposition \ref{frange} then gives
\[
\sup_{t\le(1/256)\phi(2^{n+1})-T_{-1}}r_t\le2^n
\]
for all $n\ge1$. This implies that
\[
\sup_{t\le(1/256)\phi(R)- T_{-1}}r_t\le2R
\]
for all $R\ge0$. Denote by $\psi$ the unique inverse function of
$\phi
$. Letting $R= \psi(256(T+T_{-1}))$ in the above inequality, we have
\[
\sup_{t\le T} r_t\le2\psi\bigl(256(T+ T_{-1})\bigr)\le512\psi(512T)
\]
for all sufficiently large $T$. This shows that $512\psi(512t)$ is an
upper rate function of Brownian motion on $M$ under the probability
${\mathbb P}
_{B_1}$. The technical point of passing from the average probability
${\mathbb P}_{B_1}$ to the pointwise probability ${\mathbb P}_x$ is
taken care of in
the proof of the our main theorem below.
\begin{theorem} \label{main}Let $M$ be a complete Riemannian manifold
and let $x\in M$. Let $B(R)$ be the geodesic ball on $M$ of radius $R$
and centered at $z$. Define
\[
\phi(R) = \int_6^R\frac{r \,dr}{{\ln}\vert B(r)\vert+ \ln\ln r}
\]
and let $\psi$ be the inverse function of $\phi$. There then exists a
constant $C$ such that $C\psi(C t)$ is an upper rate function of
Brownian motion $X$ on $M$, that is,
\[
{\mathbb P}_x\{ d(X_t, x)\le C\psi(C t)\mbox{ for all sufficiently large
}t\}= 1.
\]
\end{theorem}
\begin{pf} Let
\[
H = \{ d(X_t, X_0)\le C\psi(C t)\mbox{ for all sufficiently large
}t\}.
\]
We have shown that ${\mathbb P}_{B_1}(H) = 1$. This shows that $C\psi
(C t)$ is
an upper rate function for Brownian motion on $M$ starting from the
uniform distribution on the geodesic ball $B_1$. Passing to a single
starting point is easy. Let
\[
h(z) = {\mathbb P}_z(H).
\]
Let $\theta_t\dvtx {\mathscr P}(M)\rightarrow{\mathscr P}(M)$ be the
shift operator
defined by
\[
(\theta_t\omega)(s) = \omega(s+t).
\]
By the definition of the event $H$, it is clear that for any stopping
time $\tau$, we have $\omega\in H$ if and only if $\theta_\tau
\omega\in
H$; in other words, $I_H\circ\theta_\tau= I_H$. It follows that
$h(z) = {\mathbb P}_z(H) = {\mathbb E}_z I_H$
is a harmonic function on $M$. On the other hand, we have $0\le h\le1$ and
\[
\frac1{\vert B_1\vert}\int_{B_1} h(z) \,dz = 1.
\]
By the maximum principle for harmonic functions, we see that $h$ must
be identically equal to 1.
\end{pf}

The following special cases have all appeared in the literature (see
the references cited in Section \ref{sec1}). They now follow from our main
Theorem \ref{main} and all are now valid without any
geometric restrictions.
\begin{corollary} \label{oldrates} Let $M$ be a complete Riemannian
manifold. Under the following volume growth conditions, $\psi$ is an
upper rate function for Brownian motion on $M$:

(1) $\vert B(r)\vert\le C r^D$ and $\psi(t) = C_1\sqrt{t\ln t}$;

(2) $\vert B(r)\vert\le e^{C r^\alpha}$ $(0<\alpha< 2)$ and $\psi
(t) = C_1 t^{1/(2-\alpha)}$;

(3) $\vert B(r)\vert\le e^{Cr^2}$ and $\psi(t)= C_1 \exp
(C_1t^2\ln t )$;

(4) $\vert B(r)\vert\le e^{Cr^2\ln r}$ and $\psi(t)= \exp
(\exp
(C_1t ) )$.
\end{corollary}
\begin{pf} These upper rate functions follow directly from the main
theorem. Since the volume grows faster than the additional term $\ln
\ln
r$ in the function $\phi$, these rate functions are the same as if the
additional term were not there.
\end{pf}

Riemannian manifolds with slow volume growth are interesting test cases
for our main result Theorem \ref{main}. Although, in general, slow
volume growth corresponds to slow upper rate functions, our result will
not yield upper rate functions better than $\sqrt{t\ln\ln t}$ once
$\vert B(r)\vert\le(\ln r)^{\gamma}$ or some $\gamma>0$ (see Remark
\ref{nobetter}). It also explains the condition $v(r)\ge(\ln r)^\gamma$ in
the following result, due to Grigor'yan \cite{Grigoryan2}.
\begin{corollary} \label{newrate} Let $M$ be a complete manifold such
that $\vert B(r)\vert\le v(r)$ for an increasing function $v(r)\ge
(\ln
r)^\gamma$ with some $\gamma> 0$. Define $R(t)$ by
\[
\frac{R(t)^2}{\ln v(R(t))} = t.
\]
Then $CR(Ct)$ is an upper rate function for Brownian motion on $M$. In
particular, if $M$ has finite volume, then $C\sqrt{t\ln\ln t}$ is an
upper rate function.
\end{corollary}
\begin{pf} With the lower bound for $v(r)$, we have
\[
\phi(R) = \int_6^R\frac{r \,dr}{{\ln}\vert B(r)\vert+ \ln\ln r}\ge
C_1\int_6^R \frac{r \,dr}{\ln v(r)} \ge\frac{C_2R^2}{\ln v(R)}.
\]
Therefore, the inverse function $\psi(t) \le C_3R(C_3t)$. By Theorem
\ref{main}, $C R(Ct)$ is an upper rate function for some $C$.
\end{pf}
\begin{remark} In all of the concrete cases we have mentioned thus far,
upper rate functions are determined up to multiplicative constants. The
question naturally arises as to whether we could have been more careful
in our computations so as to recover the best constants in some cases,
for instance, $\psi (t) = \sqrt{(2+\varepsilon)t\ln\ln t}$ for the
standard one-dimensional Brownian motion.
Such a precise upper rate function
is impossible without
further geometric assumptions other than the volume growth. This can be
explained by means of manifolds with power volume growth $\vert
B(r)\vert\le Cr^D$. According to Corollary \ref{oldrates}(1), the
corresponding rate function is $\psi(t) = C_1\sqrt{t\ln t}$. By
comparison with a Euclidean Brownian motion, we would expect a double
logarithm instead of a single one. However, there are known examples
showing that the above rate function with a single logarithm is indeed
sharp up to a multiplicative constant (see Grigor'yan and Kelbert \cite
{GrigoryanK2}). This is the reason why we have been somewhat cavalier
about multiplicative constants in our proofs. It should be pointed out
that these constants, denoted by $C$ with or without subscripts, are
universal; they do not depend on the manifold $M$ (not even on its dimension).
\end{remark}

\section*{Acknowledgments}
The first author would like to acknowledge the hospitality and
financial support provided during his visit to the Institute of
Applied Mathematics of the Academy of Mathematics and Systems Science
at the Chinese Academy of Sciences in the summer of 2009, during which
part of this research was conducted.

%

%
\printaddresses

\end{document}